\documentclass[12pt,oneside]{amsart}
\usepackage{amsaddr}
\usepackage{graphicx}
\usepackage{enumerate}
\usepackage[round]{natbib}
\usepackage{tabularx}
\usepackage[english]{babel}
\usepackage{dsfont}
\usepackage{amssymb}
\usepackage{amsmath}
\usepackage{amsfonts}
\usepackage{hhline}
\usepackage{multirow}

\usepackage{color}

\numberwithin{equation}{section}
\newtheorem{theo}{Theorem}[section]
\newtheorem{ass}{Assumption}[section]
\newtheorem{lem}{Lemma}[section]

\newtheorem{prop}{Proposition}[section]
\newtheorem{rmk}{Remark}[section]
\newtheorem{defi}{Definition}[section]

\newcommand{\argmin}[1]{\underset{#1}{\operatorname{arg}\!\operatorname{min}}\;}
\newcommand{\abs}[1]{\lv #1\rv}

\newcommand{\cal}{\mathcal}
\newcommand{\cov}{{\mrm{Cov\,}}}

\newcommand{\de}{\delta}

\newcommand{\esp}{{\mathds{E}\,}}
\newcommand{\espt}{{\mathds{E}_{\tht}}}
\newcommand{\ex}{\exists\,}

\newcommand{\indi}{\mathds{I}}

\newcommand{\lv}{\left\vert}
\newcommand{\mrm}{\mathrm}

\newcommand{\norm}[1]{\left\Vert#1\right\Vert}
\newcommand{\p}{\mathds{P}}

\newcommand{\reel}{{\mathds{R}}}

\newcommand{\rv}{\right\vert}
\newcommand{\sste}{\scriptscriptstyle}
\newcommand{\tht}{\theta}

\newcommand{\tv}{\xrightarrow}

\newcommand{\ve}{\varepsilon}

\usepackage{setspace}
\title[]{Weighted M-estimators for multivariate clustered data: theory and simulation results}
\author[M. El Asri, D. Blanke, E. Gabriel]{Mohammed El Asri, Delphine Blanke, Edith Gabriel}
\address[M. El Asri, D. Blanke and E. Gabriel]{Avignon University, LMA EA2151 \\ Campus Jean-Henri Fabre, 301 rue Baruch de Spinoza,  \\ BP 21239, F-84916 AVIGNON Cedex 9.}

\email[E. El Asri]{mohammed.elasri@gmail.com}
\email[D. Blanke]{delphine.blanke@univ-avignon.fr}
\email[E. Gabriel]{edith.gabriel@univ-avignon.fr}

\begin{document}

\hyphenation{sa-tisfies po-sitive stu-died}

\begin{abstract}
We study weighted M-estimators  for $\reel^d$-valued clustered data and  give sufficient conditions for their consistency. Their asymptotic normality is established with estimation of the asymptotic covariance matrix. We address the robustness of these estimators in terms of their breakdown point. Comparison with the unweighted case is performed with some numerical studies. They highlight that optimal weights maximizing the relative efficiency  have a bad impact on the breakdown point.
\end{abstract}

\keywords{M-estimation, Clustered data, Efficiency, Breakdown point}

\subjclass[2000]{62F35, 62G35}

\thanks{ANR SYSCOMM, ROLSES `Robust and Optimal Locations for Sustainable Environment and Systems'.}
\maketitle


\section{Introduction}
M-estimators were first introduced by  \citet{Hu64}  as robust estimators of location and gave rise to a substantial literature. For results on their asymptotic behavior and robustness (with the study of the influence function and the breakdown point), two reference books are those of \citet{Hu81} and \citet{Ha86}. Concerning more recent works, we may refer to \citet{RG12} for  a nice  introductory presentation of robust statistics, and  to the book of \citet{Va00} for results in the multivariate case.

Most of references stand for independent and identically distributed variables. However clustered and hierarchical data
frequently arise in applications. Typically the facility location problem is an important research topic
in spatial data analysis for the geographic location of some economic activity. In this field, recent studies perform spatial modelling with
clustered data \citep[see e.g.][and references therein]{LG08,JS14}. For the multivariate one-sample location problem with clustered data,  \citet{NLO06} study the spatial median for robust estimation.  They show that the intra-cluster correlation has an impact on the asymptotic covariance matrix.  The weighted spatial median, introduced in their pioneer paper of \citeyear{NLO07},  has a better efficiency with respect to its unweighted version, especially when clusters' sizes are heterogeneous  or in the presence of strong intra-cluster correlation.   \citet{EA13} introduces a class of weighted M-estimators which generalizes their work  to a boarder class of estimators. Weights are assigned to the objective function that defines M-estimators, the aim is to adapt them to the clustered framework in order to increase their efficiency and/or robustness.

We begin with the consistency of weighted M-estimators and establish their asymptotic normality. Also, we provide consistent estimators of the asymptotic variance-covariance matrix. It allows us to derive numerically optimal weights that improve the efficiency of estimation. Finally, from a weight-based formulation of the breakdown point, we illustrate how these optimal weights may lead to its deterioration.
\section{The framework}

Clusters appear in various domains:  one may refer  to  \citet{LD02} for  current approaches in statistical inference. The framework consists in aggregated data with a homogeneous internal structure. Here, we consider $n$ independent clusters, $X_1,\dots,X_n$, where $X_i$ is made up of $m_i$ $\reel^d$-valued random variables $X_{ij}$, $j=1,\dots,m_i$, $m_i \geq 1$. Each $X_{ij}$ has the same distribution $\p_{\theta}$, $\tht\in \Theta \subset\reel^d$, where $\Theta$ is a convex and bounded set with non empty interior. The latter condition is strong, but it allows to derive results in the multivariate case for a general class of weighted M-estimators. It can be dropped when one considers some specific estimators, as for example the spatial median. For each cluster $i$, $i=1,\dotsc,n$, we also make the assumption that $(X_{i1},X_{i2})\stackrel{d}{=}(X_{ik},X_{ik'})$  for all  $k,k'=1,..,m_i$ with $k\neq k'$. This condition implies that the intra-correlation, say $r_i$, of a given cluster $X_i$ is the same for all pairs of variables inside but
 may vary from one cluster to another.  The total number of observations is denoted by $N_n:=\sum^n_{i=1}m_i$ and we suppose that  $\lim_{n \to \infty} \frac{N_n}{n}=\ell$, with  $\ell \in [1,\infty [$. Finally, $\tht$  satisfies
$$\theta= \argmin{a\in\Theta} M(a):= \argmin{a\in\Theta}\espt(\rho(X_{11},a))=\argmin{a\in\Theta} \int{\rho(x,a)dP_{\theta}(x)},$$
where, for all $a \in \Theta$,  $a \mapsto \rho(x,a)$ is a measurable function in $x$ from $\mathds{R}^d\times\mathds{R}^d$ to $\mathds{R}$.

\begin{defi}
\label{d21}
The weighted M-estimator associated with the function $\rho$  is defined by:
$$\widehat{\theta}_n^w= \argmin{a\in\Theta}  M_n^w(a) \text{ with } M_n^w(a)=\frac{1}{N_n}\sum^n_{i=1}\sum^{m_i}_{j=1} w_{ij}{ \rho(X_{ij},a)}.$$
The values $w_{ij}$, $i=1,\dotsc,n$, $j=1,\dotsc, m_i$ represent the (non random) positive weights.
\end{defi}

If $a\mapsto \rho(\cdot,a)$ is differentiable for $a=(a_1,\dotsc,a_d)^{\sste T}$ in a neighborhood of $\theta$, then the vector of partial derivatives $\psi=(\frac{\partial{\rho} }{\partial{a_1}},\dotsc,\frac{\partial{\rho} }{\partial{a_d}})^{\sste T}$ is such that  $\espt(\psi(X_{11},\theta))=0$. This leads to an alternative formulation of $\widehat{\theta}_n^w$.
\begin{defi}
\label{d22}
The weighted M-estimator $\widehat{\theta}_n^w$ is the value of `$a\!$' satisfying the $d$-vectorial equality: $
\sum^n_{i=1}\sum^{m_i}_{j=1} w_{ij}{\psi(X_{ij},a)}=0.
$
\end{defi}
Several choices are possible for the weights and  $w_{ij}\equiv 1$ returns the unweighted estimators. A particular interesting case corresponds to clusters with the same weight, namely  $w_{ij}\equiv w_i$, where the choice of $w_i$ depends on the considered framework. Typically, one may choose  $w_i = c/m_i$, $c>0$, to penalize biggest clusters or $w_i = c/r_i$ to reduce the intra-correlation effect.

\section{Asymptotic results}

Our asymptotic results  are derived under the following conditions (recall that $N_n=\sum_{i=1}^n m_i$ and $N_n/n \to \ell \ge 1$).
\begin{ass}[A\ref{a31}] Let us assume that:
\label{a31}
\begin{enumerate}[(a)]
\item For  all $\epsilon >0$:
$\underset{a \in \Theta: \left\|a-\theta\right\|\ge \epsilon}{\text{inf}}{\espt(\rho(X_{11},a))}>\espt(\rho(X_{11},\theta));$
\item  $\lim\limits_{n\to\infty}\sum\limits^n_{i=1}\big(\frac{\sum\limits^{m_i}_{j=1} w_{ij}}{N_n}\big) = 1$;
\item   $\sup\limits_{a\in\Theta}\espt(\rho^{1+\delta}(X_{11},a))<\infty$ for some $\delta$ with
\begin{enumerate}[(i)] \item either $\delta>0$ and    $\lim\limits_{n\to\infty} \sum\limits_{i=1}^n \big(\frac{\sum\limits_{j=1}^{m_i}  w_{ij}}{N_n}\big)^{1+\delta}= 0$;
\item or $\delta>1$,  $\sum\limits_{n\ge 1}\sum\limits_{i=1}^n \big(\frac{\sum\limits_{j=1}^{m_i}  w_{ij}}{N_n} \big)^{1+\delta}<\infty$ and   $\sum\limits_{n\ge 1}\big(\sum\limits_{i=1}^n \big(\frac{\sum\limits_{j=1}^{m_i}  w_{ij}}{N_n} \big)^{2}\big)^s<\infty$ for some $s>0$.
\end{enumerate}

\end{enumerate}
\end{ass}
The condition A\ref{a31}(a) guarantees the uniqueness of   $\theta$. It is satisfied if  $\rho$ is a strictly convex function and the support of   $\p_{\theta}$ is not concentrated on a line \citep[see][for the spatial median]{MD87}. The condition A\ref{a31}(c), linked both with function $\rho$ and distribution $\p_{\tht}$,  is standard and is  involved in the consistency of the estimator.  Conditions A\ref{a31}(b)-(c)  depend on the sizes $m_i$ and the weights $w_{ij}$. A\ref{a31}(b) is straightforward in the unweighted case ($w_{ij}\equiv 1$) whereas A\ref{a31}(c)-(i) follows easily for e.g. bounded $m_i$. For the case $w_{ij}\equiv w_i$ with $w_n\to 1$, the  Ces\`{a}ro theorem gives A\ref{a31}(b)  and again,  A\ref{a31}(c)-(i) is fulfilled for  bounded sequences $(m_i)$ and $(w_i)$.  The second condition A\ref{a31}(c)-(ii) is the most stringent one and it is involved in the almost sure convergence.  Remark that the choice $w_i = \frac{\ell}{m_i}$  fulfills all the conditions. Now, we may derive consistency of our estimators given in Definition~\ref{d21}.
\begin{theo}
\label{t31}
Suppose that conditions A\ref{a31}(a)-(c) are satisfied and that, for all $x$, the function $a\mapsto\rho(x,a)$ is $k(x)$-H\"{o}lderian:  $$\abs{\rho(x,a_1) - \rho(x,a_2)} \le k(x)\norm{a_1-a_2}^{\lambda},  \; 0<\lambda\le 1, \; a_1,a_2\in\Theta$$ with $\espt(k^{1+\delta}(X_{11}))<\infty$ and $\delta$ given in A\ref{a31}(c). Then, the condition A\ref{a31}(c)-(i) gives the $P$-consistency of $\widehat{\theta}_n^w$  while  A\ref{a31}(c)-(ii) gives its strong consistency.
\end{theo}

\begin{proof}
With Assumption \ref{a31}, the relation  $\left\|a-\tht\right\|\ge \epsilon$, for $a\in\Theta$ and some $\ve >0$,  implies that $M(a)>M(\theta)+\eta$ for some $\eta>0$. So, one gets the inclusion: $\left\{\left\|\widehat{\theta}_n-\theta\right\| \ge \epsilon\right\}\subset \left\{ M(\widehat{\theta}_n)>M(\theta)+\eta\right\}$. By this way,
$$
\left\{0\leq M(\widehat{\theta}_n^w)-M(\theta)\leq \eta\right\} \subset \left\{\left\|\widehat{\theta}_n^w-\theta\right\| < \epsilon\right\}.
$$
Moreover, the estimator  $\widehat{\theta}_n^w$ is the minimizer of $M_n^w$, hence $M_n^w(\widehat{\theta}_n^w)\leq M_n^w(\theta)$ and we get:
\begin{multline*}
0 \le   M(\widehat{\theta}_n^w)-M(\theta) = M(\widehat{\theta}_n^w)-M_n^w(\widehat{\theta}_n^w)+M_n^w(\widehat{\theta}_n^w)-M(\theta)
\\\le  M(\widehat{\theta}_n^w)-M_n^w(\widehat{\theta}_n^w)+M_n^w(\theta)-M(\theta)
\le 2 \sup_{a \in \Theta}{\left|M_n^w(a)-M(a)\right|}.
\end{multline*} Hence,   the inclusion $
\left\{\sup_{a \in \Theta} \left|M_n^w(a)-M(a)\right| <\eta\right\} \subset \{\|\widehat{\theta}_n^w-\theta\| < \epsilon\}$.  The result follows from the uniform consistency of $M_n^w$  established in the next lemma (with proof postponed to the end).
\begin{lem}
\label{l31}
Under conditions of Theorem~\ref{t31}, $M_n^w(a)-M(a)$ converges to $0$ uniformly in $a\in\Theta$, either in probability or almost surely.
\end{lem}\end{proof}
\begin{rmk}\label{rem1}
Theorem~\ref{t31} includes estimators with lipschitzian objective function as the weighted spatial median with $\rho(x,a) = \norm{x-a}$ and the weighted Huber estimator derived  from
 $\rho(x,a) = \frac{1}{2}\norm{x-a}^2$ if   $\norm{x-a} \le k$ and $\rho(x,a)= k \norm{x-a}-\frac{1}{2} k^2$ if $ \norm{x-a} > k$.\end{rmk}

We conclude this part with the asymptotic normality of our estimators. To this aim, we use the following conditions where $\dot{\psi}$ denotes the Hessian of the function  $a\mapsto \rho(x,a)$.
\begin{ass}[A\ref{a32}] Let us assume that:
\label{a32}
\begin{enumerate}[(a)]
\item  $\lim\limits_{n\to\infty}\sum\limits^n_{i=1}\big(\frac{\sum\limits^{m_i}_{j=1} w_{ij}}{N_n}\big)=1$ and $\lim\limits_{n\to\infty}{\frac{1}{N_n}\sum\limits^n_{i=1}\sum\limits^{m_i}_{j=1}{w_{ij}^2}}=c_{w}<\infty$;
\item  $\lim\limits_{n\to \infty}{\frac{1}{N_n}\sum\limits^n_{i=1}\sum\limits^{m_i}_{j=1} \sum\limits_{j'\neq j}{w_{ij}w_{ij'}C_{\tht,i}}}=C_{\tht}^w$ with $C_{\tht,i}=\espt\psi(X_{ij},\theta)\psi^{\sste T}(X_{ij'},\theta)$;
\item $\espt(\left\|\psi(X_{11},\theta)\right\|^{2+\eta})<\infty$   for some $\eta>0$ and $\lim\limits_{n\to \infty} \frac{1}{N_n^{1+\frac{\eta}{2}}}\sum\limits_{i=1}^n (\sum\limits_{j=1}^{m_i} w_{ij})^{2+\eta}= 0$;
\item  $\espt\big(\|\dot{\psi}(X_{11},\theta)\|\big)^{1+\de}<\infty$  for some $\de>0$ (with matrix-norm)  and $a\mapsto \dot{\psi}(x,a)$ is $k(x)$ - H\"{o}lderian: $\big\|\dot{\psi}(x,a_2) - \dot{\psi}(x,a_1)\big\| \le k(x)\norm{a_2-a_1}^{\lambda}$ with $\esp(k^{1+\de}(X_{11}))<\infty$ and $\lim\limits_{n\to\infty}\sum\limits_{i=1}^n \big(\frac{\sum_{j=1}^{m_i}  w_{ij}}{N_n}\big)^{1+\delta} = 0$.
\end{enumerate}
\end{ass}
If $w_{ij}\equiv w_i$, with $\lim_{n\to \infty} w_n=1$, the condition A\ref{a32}(a) is fulfilled. Moreover if the intra-cluster correlations are the same for all clusters, $C_{\tht,i}=C_{\tht}$, and if  $\lim_{n\to\infty}m_n=\ell$,  then one gets  $C_{\tht}^w=(\ell-1)C_{\tht}$: the condition A\ref{a32}(b) is also satisfied. Condition A\ref{a32}(c) is involved for the application of the Lindeberg theorem and  its last part is fulfilled with the previous choices of weights. Also,  it holds true as soon as $\frac{1}{n}\sum_{i=1}^n (\sum_{j=1}^{m_i} w_{ij})^{2+\eta}$ is asymptotically finite. The weighted spatial median is not covered by the condition A\ref{a32}-(d), but its study may be direct, as done in  \citet{NLO06}. Here our conditions are satisfied with sufficiently smooth $\psi$ functions: in particular the H\"{o}lderian condition holds as soon as each element of $\dot{\psi}$ is $k(x)$-H\"{o}lderian. By this way, we do not require the second derivative of $\psi$, which seems reasonable from a robust point of view. In particular, $L_p$-medians defined by $\widehat{\tht}_n^w = \argmin{a\in\Theta} \frac{1}{N_n} \sum_{i=1}^{n} \sum_{j=1}^{m_i} w_{ij} \norm{X_{ij} - a}^p$ fulfill these conditions with $p \ge 2$.

\begin{theo}
\label{t32}
Suppose that $\widehat{\theta}_n^w \tv[n\to\infty]{p} \tht$ and that Assumption  \ref{a32} is fulfilled. Also assume  that  $V_{\theta}:=\espt(\dot{\psi}(X_{11},\tht))$ is invertible. Then,
$$
\sqrt{N_n}(\widehat{\theta}_n^w-\theta) \xrightarrow[n\to\infty]{d} N\left(0,\Sigma_{\tht}^w\right)
$$
where $\Sigma_{\tht}^w=V_{\theta}^{-1}\left(c_{w}B_{\tht}+C_{\tht}^w\right)V_{\theta}^{-1}$ and $B_{\tht}:=\espt\psi(X_{11},\theta)\psi^{\sste T}(X_{11},\theta)$.
\end{theo}

\begin{proof}
First, for $a=(a_1,\dotsc,a_d)^{\sste T}$ in a convex open neighborhood of $\theta$, let us denote
$$
\psi(x,a)=(\frac{\partial{\rho(x,a)}}{\partial{a_1}},\dotsc,\frac{\partial{\rho(x,a)}}{\partial{a_d}})^{\sste T} =:(\psi_1(x,a),\dotsc,\psi_d(x,a))^{\sste T}.$$
We apply a Taylor expansion with integral remainder:
$$
\psi(x,a)=\psi(x,\theta)+\int_0^1 \dot{\psi}(x,\theta+s(a-\theta))^{\sste T}(a-\theta)\, \mathrm{d}s,
$$
For  $a=\widehat{\theta}_n^w\xrightarrow[n\to\infty]{p}\theta$ and $x=X_{ij}$, we get:
\begin{multline*}
\frac{1}{N_n}\sum_{i=1}^n\sum_{j=1}^{m_i} w_{ij} \psi(X_{ij},\widehat{\theta}_n^w)= \frac{1}{N_n} \sum_{i=1}^n \sum_{j=1}^{m_i} w_{ij} \psi(X_{ij},\theta)\\+\big(\frac{1}{N_n}\sum_{i=1}^n \sum_{j=1}^{m_i} w_{ij} \dot{\psi}(X_{ij},\theta)^{\sste T} + o_p(1)\big) (\widehat{\theta}_n^w-\theta)
\end{multline*}
where $o_p(1)$ comes from  the $P$-consistency of $\frac{1}{N_n} \sum_{i=1}^n\sum_{j=1}^{m_i} w_{ij}k(X_{ij})$, derived from  Theorem~\ref{ta1}(1) with conditions A\ref{a32}(a) and (d), together with  the H\"{o}lderian condition  made on  $\dot{\psi}$.

Next, if $T_n^w(a) = \frac{1}{N_n}\sum_{i=1}^n \sum_{j=1}^{m_i} w_{ij} \psi(X_{ij},a)$, then $T_n^w(\widehat{\theta}_n^w)=0$ by definition of $\widehat{\theta}_n^w$ and we obtain
\begin{equation*}
  -\sqrt{N_n}T_n^w(\theta) =\left(\dot{T}_n^w(\theta)+o_p(1)\right)\sqrt{N_n}(\widehat{\theta}_n^w-\theta).
\end{equation*}
Next, the two following technical lemmas (with proofs  postponed to the end) and the $P$-consistency of $\widehat{\tht}^w_n$ yield the result.
\begin{lem}
\label{l32}
Under the conditions A\ref{a32}(a)-(c),
$$\sqrt{N_n}T_n^w(\theta) \xrightarrow[n\to\infty]{d}N\left(0,c_{w}B_{\tht}+C_{\tht}^w\right).$$
\end{lem}
\begin{lem}
\label{l33}
If assumptions of Theorem \ref{t32} hold, we get $$\dot{T}_n^w(\theta) \xrightarrow[n\to\infty]{p}V_{\theta}.$$
\end{lem}
\end{proof}

In the independent case (obtained with $m_i\equiv 1$), the asymptotic variance is reduced to $c_wV_{\theta}^{-1}B_{\tht}V_{\theta}^{-1}$ with  $c_w$  equal to 1 for the unweighted case ($w_{ij}\equiv 1$). We may deduce that weights are of no use to improve this variance. Actually, the condition A\ref{a32}(a) gives $\frac{1}{n}\sum_{i=1}^n w_i \tv[n\to\infty]{}1$, so by the Cauchy-Schwarz inequality, we get that $c_w= \lim\limits_{n\to\infty} \frac{1}{n}\sum_{i=1}^n w_i^2$ is  always greater than one. Clustering effects appear through the term $C_{\tht}^w$, so we may  choose the weights to reduce this value in relation to the unweighted case.  Finally in the case of equal intra-correlation ($C_{\tht,i}\equiv C_{\tht}$), it should be noticed that  the choice minimizing the variance depends neither on $\p_{\tht}$ nor on the objective function $\rho$.

\section{Relative efficiency}

\subsection{Estimation of the asymptotic variance}
We study the relative efficiency of weighted M-estimators   toward their  unweighted versions. First, recall that the unweighted case is obtained with $w_{ij}\equiv 1$ for which $\Sigma_{\tht}^w:=\Sigma_{\tht} = V_{\tht}^{-1} (B_{\tht}+C_{\tht}^m)V_{\tht}^{-1}$, as $c_w=1$ by condition A\ref{a32}-(a) and $C_{\tht}^w:=C_{\tht}^m = \lim_{n\to\infty} \frac{1}{N_n}\sum_{i=1}^n m_i(m_i-1)C_{\tht,i}$. Variances are compared by using the relative efficiency index defined by
\begin{equation}
\label{e35}
E_f=\bigg[\frac{\text{det}(V_{\theta}^{-1}\left(B_{\tht}+C_{\tht}^m\right)V_{\theta}^{-1})} {\text{det}(V_{\theta}^{-1}\left(c_{w}B_{\tht}+C_{\tht}^w\right)V_{\theta}^{-1})}\bigg]^{\frac{1}{d}} = \bigg[\frac{\text{det}(\Sigma_{\tht})} {\text{det}(\Sigma_{\tht}^w)}\bigg]^{\frac{1}{d}}.
\end{equation}
First we propose estimators for these variances and study their  consistency under assumptions of Theorem~\ref{t32}.  To estimate $c_wB_{\tht}$, we denote $\widehat{B}_n^w(a)$ the matrix defined by:
$$
\widehat{B}_n^w(a) =\frac{1}{N_n}  \sum_{i=1}^n  \sum_{j=1}^{m_i} w_{ij}^2  \psi(X_{ij},a)\psi^{\sste T}(X_{ij},a).
$$
Similarly, estimators of  $C_{\tht}^w$ and $V_{\theta}$ are derived from the functionals
$$
\widehat{C}_n^w(a)=\frac{1}{N_n}  \sum_{i=1}^n   \sum_{j=1}^{m_i} \sum_{j'\neq j }^{m_i} w_{ij}w_{ij'}\psi(X_{ij'},a)\psi^{\sste T}(X_{ij},a)$$
and $\widehat{V}_n(a)=\frac{1}{N_n}  \sum_{i=1}^n  \sum_{j=1}^{m_i}  w_{ij}\dot{\psi}^{\sste T}(X_{ij},a)$. Consistency of these means is derived in the following proposition with proof postponed to the Appendix.
\begin{prop}\label{p41}
Suppose that assumptions of  Theorem \ref{t32}  are fulfilled.
\begin{enumerate}[1)]
\item If  $\lim\limits_{n\to\infty}\sum\limits_{i= 1}^n \big(\frac{\sum_{j=1}^{m_i} w_{ij}^2}{N_n}\big)^{1+\frac{\eta}{2}}=0$, then  $\widehat{B}_n^w(a) \tv[n\to\infty]{p}c_wB_{a}$,  $a\in\Theta$.
\item If  $\lim\limits_{n\to\infty}\sum\limits_{i= 1}^n \big(\frac{\sum_{j=1}^{m_i} \sum_{j'=1,j'\not=j}^{m_i}w_{ij}w_{ij'}}{N_n}\big)^{1+\frac{\eta}{2}}=0$, then $\widehat{C}_n^w(a) \tv[n\to\infty]{p}C_{a}^w$,  $a\in \Theta$.
\item  $\widehat{V}_n(a) \tv[n\to\infty]{p} V_{a}$, $a\in \Theta$.
\end{enumerate}
\end{prop}
Note that these additional conditions are again satisfied for the weights discussed in the previous section. Next, under some additional regularity assumptions on  $\psi$ and its derivatives, the previous results turn to be uniform on $a\in \Theta$. This allows us to construct consistent estimators of the asymptotic variance-covariance matrix given in Theorem~\ref{t32} by plugging $\widehat{\theta}^w_n$ in $\widehat{B}_n^w(a)$,  $\widehat{C}_n^w(a)$, $\widehat{V}_n(a) $. We get the following result.
\begin{theo}\label{t41} Suppose that assumptions of Proposition~\ref{p41} are fulfilled. If moreover,
\begin{enumerate}[1)]
\item $\sup\limits_{a\in\Theta} \espt\big\|\psi(X_{11},a)\big\|^{2+\eta}<\infty$ and $$\big\| \psi(x,a_2)\psi^{\sste T}(x,a_2) - \psi(x,a_1)\psi^{\sste T}(x,a_1) \big\| \le k(x) \big\|a_2-a_1\big\|^{\lambda},\;\;\lambda\in ]0,1] $$ with $\espt (k^{1+\frac{\eta}{2}} (X_{11}))< \infty$, then $\widehat{B}^w_n(\widehat{\tht}_n^w) \tv[n\to\infty]{p} c_w B_{\tht}$;
\item $\sup\limits_{a\in\Theta} \espt\big\|\psi(X_{11},a)\big\|^{2+\eta}<\infty$ and $$\big\| \psi(x,a_2)\psi^{\sste T}(y,a_2) - \psi(x,a_1)\psi^{\sste T}(y,a_1) \big\| \le k(x,y) \big\|a_2-a_1\big\|^{\lambda},\;\;\lambda\in ]0,1] $$ with $\espt (k^{1+\frac{\eta}{2}}(X_{i1},X_{i2})) < \infty$ for all $i=1,\dotsc,n$ and $$\lim\limits_{n\to\infty} \frac{1}{N_n} \sum\limits_{i=1}^{n}\sum\limits_{j=1}^{m_i}\sum\limits_{j'\not=j} w_{ij}w_{ij'} \espt(k(X_{ij},X_{ij'})) = k_w<\infty,$$ then $\widehat{C}^w_n(\widehat{\tht}_n^w) \tv[n\to\infty]{p} C_{\tht}$;
\item  $\sup\limits_{a\in\Theta} \espt\big\|\dot{\psi}(X_{11},a)\big\|^{1+\delta}<\infty$ then  $\widehat{V}_n(\widehat{\tht}_n^w) \tv[n\to\infty]{p}
    V_{\tht}$.
\end{enumerate}
\end{theo}

\begin{proof}

1) Concerning the term  $\widehat{B}^w_n(\widehat{\tht}_n^w)$, we may write
\begin{align*}
\norm{\widehat{B}^w_n(\widehat{\tht}_n^w) - c_wB_{\tht}} & \le \norm{\widehat{B}^w_n(\widehat{\tht}_n^w) - c_wB(\widehat{\tht}_n^w)}+c_w\norm{B(\widehat{\tht}_n^w) - B_{\tht}}
\\
&\le \sup\limits_{a\in\Theta} \norm{\widehat{B}^w_n(a) -c_w B(a )} + c_w\norm{B(\widehat{\tht}_n^w) - B_{\tht}}.
\end{align*}
The H\"{o}lderian condition on $a\mapsto \psi(x,a)\psi^{\sste T}(x,a)$ gives the continuity of the function $a\mapsto \esp(\psi(X_{11},a)\psi^{\sste T}(X_{11},a))$; so the second norm tends to 0 by the $P$-consistency of $\widehat{\tht}_n^w$ and Slutsky's lemma. It remains  to establish the uniform $P$-consistency of $\widehat{B}^w_n(a)$. This result is derived similarly to the proof of Lemma~\ref{l31}: we  use again the H\"{o}lderian condition on $\psi\psi^{\sste T}$ by noting that  obtained bounds in the proof of Proposition~\ref{p41} are uniform on $a$ since $\sup\limits_{a\in \Theta} \espt\big\|\psi(X_{11},a)\big\|^{2+\eta}<\infty$.

2) The two other terms are studied similarly with the help of:
\begin{align*}
\norm{\widehat{C}_n^w(\widehat{\tht}_n^w) - C^w_{\tht}} & \le \sup\limits_{a\in\Theta} \norm{\widehat{C}_n^w(a) -C^w(a )} + \norm{C(\widehat{\tht}_n^w) - C^w_{\tht}} \\
\norm{\widehat{V}_n(\widehat{\tht}_n^w) - V_{\tht}} & \le \sup\limits_{a\in\Theta} \norm{\widehat{V}_n(a) -V(a )} + \norm{V(\widehat{\tht}_n^w) - V_{\tht}}.
\end{align*}

\end{proof}
\subsection{Numerical results}
\label{subsec:efficiency}
Different simulations and analyses were performed with \emph{R} software \citep{R14}. We use an empirical version of the efficiency $E_f$ defined in \eqref{e35} with $\Sigma_{\tht}$ and $\Sigma_{\tht}^w$ respectively replaced by $\widehat{\Sigma}_{\tht}$ and $\widehat{\Sigma}_{\tht}^w$. Note that here, we use the true value $\tht$ in order to compare the asymptotic efficiencies in an ideal setting. For 1000 replications of each case, we compute four configurations of 100 multivariate random variables subdivided in 10 clusters:
\begin{enumerate}
\item C1: 9 clusters of size 4  and 1 cluster  of size 64;
\item  C2:  5 clusters of size 4 and 5 clusters of size 16;
\item  C3: 2 clusters of size 4,  1 cluster of size  8, and 7 clusters of size 12;
\item  C4: 10 clusters of all sizes  from 5 to 15 except 10.
\end{enumerate}
\begin{table}[t!]
  \centering \caption{Optimal weights for the weighted spatial median, the weighted Huber estimator and the $L_p$ estimator ($p=3$) obtained for a centered bivariate Gaussian distribution with  $r=0.2$ and $r=0.8$}
\renewcommand{\arraystretch}{1.2}
\scalebox{0.75}{\begin{tabular}{|l||r|lr||r|lr||r|lr||r|lr|}
\hline
\multirow{2}{*}{}
  &  \multicolumn{3}{|c||}{C1}    &  \multicolumn{3}{c||}{C2}      &  \multicolumn{3}{c||}{C3}    &  \multicolumn{3}{c|}{C4}    \\
\cline{2-13}
& $m_i$  & $r=0.2$ & $r=0.8$ &$m_i$ & $r=0.2$ & $r=0.8$ &  $m_i$ &  $r=0.2$ & $r=0.8$ & $m_i$  &  $r=0.2$ & $r=0.8$
 \\  \hhline{|=============|}
 \parbox[t]{2mm}{\multirow{10}{*}{\rotatebox[origin=c]{90}{Weighted spatial median}}} &
     4     & 2.251 & 2.487& 4     &  1.818& 2.412 & 4     & 1.712 & 2.382& 5     & 1.519 &1.939\\
&    4     & 2.324 & 2.520 & 4 &1.877   & 2.444 & 4     &1.765   &  2.414& 6     &1.420   & 1.645\\
&    4     & 2.232 &2.464 & 4     & 1.803 &  2.389  & 8     & 1.180  &1.243& 7     &  1.262   & 1.407\\
&    4     &  2.203 & 2.460 & 4     &1.780 & 2.386 & 12    & 0.915   &  0.837& 8     &1.171 &  1.241\\
&    4     & 2.252 &2.473& 4     & 1.819 &  2.398 & 12    & 0.911  &  0.846 & 9     & 1.093  &  1.114 \\
&    4     &  2.291 & 2.486 & 16    &0.813 & 0.654& 12    &0.926 & 0.851 & 11    &0.989 & 0.920\\
&    4     & 2.227  & 2.482 & 16    & 0.789 & 0.656 & 12    & 0.907   &  0.855 & 12    & 0.904  & 0.852\\
&    4     & 2.192  &  2.454 & 16    & 0.769 & 0.638  & 12    &0.900   & 0.831  & 13    &0.841  &0.767 \\
&    4     & 2.281  & 2.498 & 16    &0.821 & 0.658  & 12    &  0.929  & 0.859  & 14    & 0.840   &  0.738 \\
&    64    & 0.297   & 0.167 & 16    & 0.784  &0.636 & 12    &0.901  & 0.827 & 15    &0.763   & 0.665 \\
   \hhline{|=============|}
  \parbox[t]{2mm}{\multirow{10}{*}{\rotatebox[origin=c]{90}{Weighted-Huber}}} &
    4     & 2.268 & 2.464 & 4     & 1.875 &2.416 & 4     & 1.760 &2.393 & 5     & 1.567 &1.941\\
&    4     & 2.343 & 2.520 & 4 & 1.937  & 2.470 & 4     & 1.816   & 2.447 & 6     & 1.439  & 1.659\\
&    4     & 2.275 & 2.468 & 4     & 1.881 &   2.419 & 8     & 1.199  & 1.239 & 7     & 1.301   & 1.411\\
&    4     & 2.216 & 2.451 & 4     & 1.832 & 2.402 & 12    & 0.902  & 0.830 & 8     & 1.171  & 1.227\\
&    4     & 2.259 & 2.464 & 4     & 1.868 & 2.415 & 12    & 0.894  & 0.835 & 9     & 1.089  & 1.104 \\
&    4     & 2.338 & 2.522 & 16    & 0.797 & 0.650 & 12    & 0.922 & 0.855& 11    & 0.978 & 0.928\\
&    4     & 2.289   & 2.506 & 16    & 0.776  & 0.647 & 12    & 0.910  & 0.850 & 12    & 0.902  & 0.848\\
&    4     & 2.242  & 2.441 & 16    & 0.753 & 0.631 & 12    & 0.879   & 0.829  & 13    & 0.829  & 0.765 \\
&    4     & 2.352  & 2.563 & 16    & 0.810  & 0.662 & 12    & 0.940  & 0.869 & 14    & 0.833   & 0.746 \\
&    64    & 0.276  & 0.163 & 16    & 0.767  & 0.629 & 12    & 0.896  & 0.826 & 15    & 0.749  & 0.663\\
    \hline
 \hhline{|=============|}
  \parbox[t]{2mm}{\multirow{10}{*}{\rotatebox[origin=c]{90}{Weighted $L_3$ estimator}}} &
  4     & 2.141 &  2.370& 4     & 1.756 &  2.320 & 4     & 1.654 & 2.299 & 5     & 1.485 & 1.862\\
  &      4     & 2.314   & 2.632 & 4     & 1.897  & 2.577 & 4     & 1.787  & 2.554 & 6     & 1.478 & 1.748  \\
 &   4     & 2.331  & 2.538& 4     & 1.911   & 2.485 & 8     & 1.209   & 1.274& 7     & 1.292   & 1.442\\
 &   4     & 2.311   & 2.513& 4     & 1.894  & 2.460 & 12    & 0.900  & 0.849 & 8     & 1.192   & 1.256 \\
 &   4     & 2.145    & 2.273 & 4     & 1.759  & 2.226 & 12    & 0.861  & 0.772 & 9     & 1.036  & 1.018 \\
  &  4     & 2.241  & 2.424 & 16    & 0.777   & 0.629 & 12    & 0.923  & 0.827 & 11    & 0.964   & 0.897 \\
  &  4     & 2.271    & 2.457 & 16    & 0.786  & 0.635 & 12    & 0.911  & 0.833 & 12    & 0.902  & 0.831 \\
  &  4     & 2.095   & 2.261 & 16    & 0.734   & 0.587 & 12    & 0.854   & 0.770 & 13    & 0.788   & 0.711 \\
  &  4     & 2.550   & 2.808 & 16    & 0.862   & 0.723 & 12    & 1.010  & 0.950 & 14    & 0.889   & 0.814  \\
  &  64    & 0.287  & 0.170 & 16    & 0.787   & 0.659 & 12    & 0.922  & 0.866  & 15    & 0.779  & 0.694 \\
  \hline
\end{tabular}}
\label{tab1}\end{table}
We consider zero mean bivariate models  ($\tht=(0,0)^{\sste T}$) with independent components with Gaussian or Student distribution, the latter with $\nu$ degrees of freedom ($\nu\in\{1,3\}$). We  give the same intra-cluster correlation $r_i\equiv r\in]0,1[$ to  all the clusters and weights are taken to be the same within a given cluster ($w_{ij}\equiv w_i$). Consequently, for all $j\not=j'$: $ \cov(X_{ij},X_{ij'})=r \begin{pmatrix}
 1&0\\ 0&1 \end{pmatrix}$ and, $\cov(X_{ij}, X_{ij})=\begin{pmatrix}  1&0\\ 0&1 \end{pmatrix}$ with $i=1,\dotsc,10$. The M-estimators considered in this study are: the spatial median, the empirical mean, the Huber estimator (defined in Remark~\ref{rem1}), the $L_p$-medians with $p=3,4,5$. Finally, the results are derived with $r=0.2$ and $r = 0.8$.


First we optimize the weights to minimize $\text{det}(\widehat{\Sigma}_{\tht}^w)$ with \emph{Matlab} \citeyear{M09}. Table \ref{tab1} shows that optimal weights have the same order of magnitude for clusters of the same size, whatever the estimator  and the correlation $r$. These results are  in agreement with our theoretical ones. Moreover, we observe that weights decrease as the clusters' sizes increase: as expected, big groups of aggregated data are penalized.



\begin{table}[t!]
  \centering
  \caption{Gaussian distribution:  Relative efficiencies w.r.t. to the (weighted) mean  and $\widehat{E}_f=(\frac{\text{det}(\widehat{\Sigma}_{\tht})}{\text{det}(\widehat{\Sigma}_{\tht}^w)})^{\frac{1}{2}}$}
 \renewcommand{\arraystretch}{1.2}
\scalebox{0.75}{   \begin{tabular}{|l||c|c|c||c|c|c|}
\cline{2-7}
 \multicolumn{1}{c|}{}&  \multicolumn{1}{c|}{$\mathbf{(\frac{\text{\bf det}(\widehat{\Sigma}_{\tht}^w) }{\text{\bf det}(\widehat{\Sigma}_{\text{\bf opt}}^w)})^{\frac{1}{2}}}$} & \multicolumn{1}{c|}{$\mathbf{(\frac{\text{\bf det}(\widehat{\Sigma}_{\tht})}{\text{\bf det}(\widehat{\Sigma}_{\text{\bf opt}})})^{\frac{1}{2}}}$} &  \multicolumn{1}{c|}{$\mathbf{\widehat{E}_f}$}  &  \multicolumn{1}{c|}{$\mathbf{(\frac{\text{\bf det}(\widehat{\Sigma}_{\tht}^w) }{\text{\bf det}(\widehat{\Sigma}_{\text{\bf opt}}^w)})^{\frac{1}{2}}}$} & \multicolumn{1}{c|}{$\mathbf{(\frac{\text{\bf det}(\widehat{\Sigma}_{\tht})}{\text{\bf det}(\widehat{\Sigma}_{\text{\bf opt}})})^{\frac{1}{2}}}$} &  \multicolumn{1}{c|}{$\mathbf{\widehat{E}_f}$}\\  \hhline{~|======|}
\multicolumn{1}{l|}{}&\multicolumn{6}{c|}{Configuration C1 (left: $r=0.2$, right: $r=0.8$)}\\\hline
    \textbf{Median} & 1.167& 1.068&2.293 &1.243& 1.208& 3.923\\\hline
    \textbf{Mean} & 1 & 1 & 2.507 & 1& 1 & 4.036 \\\hline
    \textbf{Huber} &1.036 & 1.017 &2.461 & 1.069  & 1.064 & 4.019\\\hline
    \textbf{Lp-Median; $p=3$} & 1.06 & 1.008  &  2.385& 1.063& 1.044 & 3.964\\\hline
    \textbf{Lp-Median; $p=4$} & 1.257 & 1.052 &2.099 & 1.275 & 1.216 & 3.849\\\hline
    \textbf{Lp-Median; $p=5$} & 1.614 & 1.133 & 1.76&1.629  & 1.488 & 3.686\\ \hhline{|=======|}
\multicolumn{1}{|l||}{} &      \multicolumn{6}{c|}{Configuration C2 (left: $r=0.2$, right: $r=0.8$)}\\ \hline
   \textbf{Median} &1.078 & 1.064 & 1.128  & 1.066 & 1.142 &1.316\\\hline
    \textbf{Mean} & 1 &1 & 1.142& 1 & 1 &  1.332\\\hline
    \textbf{Huber} & 1.029 &1.022 & 1.134 & 1.048 & 1.046 & 1.33\\\hline
    \textbf{Lp-Median; $p=3$} & 1.059 & 1.044&1.126 & 1.103&1.095 &1.323 \\\hline
    \textbf{Lp-Median; $p=4$} & 1.217 & 1.167 & 1.095 & 1.4  &1.38 & 1.313 \\\hline
    \textbf{Lp-Median; $p=5$} & 1.51& 1.393 & 1.054 & 1.986 &1.942 & 1.303 \\ \hhline{|=======|}
 \multicolumn{1}{|l||}{} &    \multicolumn{6}{c|}{Configuration C3 (left: $r=0.2$, right: $r=0.8$)}\\\hline
   \textbf{Median} &1.112  & 1.103 &  1.032 &1.129 & 1.125&  1.091\\      \hline
    \textbf{Mean} & 1 & 1 &  1.040 &1 & 1 & 1.095 \\       \hline
   \textbf{Huber} & 1.028 & 1.026 & 1.038 & 1.048& 1.048 & 1.095\\      \hline
    \textbf{Lp-Median; $p=3$} & 1.059 & 1.052 & 1.034 &1.100  & 1.095 & 1.09\\       \hline
   \textbf{Lp-Median; $p=4$} &1.211& 1.193& 1.025& 1.385 & 1.375 & 1.087  \\      \hline
    \textbf{Lp-Median; $p=5$} &1.486 &1.447 & 1.013  &1.931 & 1.911 & 1.084\\    \hhline{|=======|}
  \multicolumn{1}{|l||}{} &    \multicolumn{6}{c|}{Configuration C4 (left: $r=0.2$, right: $r=0.8$)}\\\hline
       \textbf{Median} & 1.105 & 1.096 & 1.043 & 1.130  & 1.126& 1.100 \\  \hline
     \textbf{Mean} & 1 &1 & 1.052 & 1 & 1 & 1.103  \\    \hline
   \textbf{Huber} & 1.028 & 1.025  & 1.049& 1.048&1.048 & 1.103\\   \hline
    \textbf{Lp-Median; $p=3$} & 1.059 & 1.051  & 1.045& 1.095&1.109& 1.098\\    \hline
   \textbf{Lp-Median; $p=4$} & 1.210 &1.190& 1.035 &1.365&1.356 &  1.095 \\    \hline
   \textbf{Lp-Median; $p=5$} & 1.477& 1.435 &  1.022 &1.87 & 1.852 & 1.092\\    \hline
    \end{tabular}}%
  \label{tab2}%
\end{table}%

\begin{table}[h!]
  \centering
  \caption{Student distribution with 3 df:  Relative efficiencies w.r.t. to the (weighted) median and $\widehat{E}_f=(\frac{\text{det}(\widehat{\Sigma}_{\tht})}{\text{det}(\widehat{\Sigma}_{\tht}^w)})^{\frac{1}{2}}$}
    \renewcommand{\arraystretch}{1.2}
\scalebox{0.82}{   \begin{tabular}{|l||c|c|c||c|c|c|}
\cline{2-7}
 \multicolumn{1}{c|}{}&  \multicolumn{1}{c|}{$\mathbf{(\frac{\text{\bf det}(\widehat{\Sigma}_{\tht}^w) }{\text{\bf det}(\widehat{\Sigma}_{\text{\bf opt}}^w)})^{\frac{1}{2}}}$} & \multicolumn{1}{c|}{$\mathbf{(\frac{\text{\bf det}(\widehat{\Sigma}_{\tht})}{\text{\bf det}(\widehat{\Sigma}_{\text{\bf opt}})})^{\frac{1}{2}}}$} &  \multicolumn{1}{c|}{$\mathbf{\widehat{E}_f}$}  &  \multicolumn{1}{c|}{$\mathbf{(\frac{\text{\bf det}(\widehat{\Sigma}_{\tht}^w) }{\text{\bf det}(\widehat{\Sigma}_{\text{\bf opt}}^w)})^{\frac{1}{2}}}$} & \multicolumn{1}{c|}{$\mathbf{(\frac{\text{\bf det}(\widehat{\Sigma}_{\tht})}{\text{\bf det}(\widehat{\Sigma}_{\text{\bf opt}})})^{\frac{1}{2}}}$} &  \multicolumn{1}{c|}{$\mathbf{\widehat{E}_f}$}\\  \hhline{~|======|}
\multicolumn{1}{l|}{}&\multicolumn{6}{c|}{Configuration C1 (left: $r=0.2$, right: $r=0.8$)}\\\hline
    \textbf{Median} & 1 & 1 & 2.273 & 1 & 1& 3.894\\\hline
    \textbf{Mean} &  1.947 & 2.157 & 2.518 & 1.984& 1.995& 3.917\\\hline
    \textbf{Huber} & 1.123 & 1.225 & 2.478&1.111  &1.145 & 4.014 \\  \hline
     \textbf{Lp-Median; $p=3$} & 2.01E+3&1.86E+3 &2.096 & 2.40E+3&2.35E+3& 3.814
\\    \hline
   \textbf{Lp-Median; $p=4$} &  7.14E+4 & 3.80E+4 & 1.209& 5.72E+4&4.26E+4 &  2.904 \\    \hline
   \textbf{Lp-Median; $p=5$} &  1.54E+6 & 6.91E+5    & 1.020   &  6.28E+5  &3.08E+5  &1.911 \\    \hhline{|=======|}
\multicolumn{1}{|l||}{} &  \multicolumn{6}{c|}{Configuration C2 (left: $r=0.2$, right: $r=0.8$)}\\ \hline
    \textbf{Median} & 1& 1&1.111 & 1& 1 &  1.313 \\\hline
    \textbf{Mean} &2.078 &2.106 &  1.127&2.059 &2.051 &1.307\\\hline
    \textbf{Huber} & 1.152&1.174 & 1.133& 1.152 &1.166 &1.328
    \\  \hline
     \textbf{Lp-Median; $p=3$} &2.22E+3&2.18E+3 &1.100& 2.10E+3 &2.07 E+3  &1.299 \\    \hline
   \textbf{Lp-Median; $p=4$} &  7.2ZE+4 &6.63E+4& 1.022 & 3.88e+4 & 3.51E+4  & 1.188 \\    \hline
   \textbf{Lp-Median; $p=5$} & 1.43E+6 &1.300E+6    &1.007 & 4.11E+5  & 3.41E+5 & 1.089 \\    \hhline{|=======|}
\multicolumn{1}{|l||}{} &  \multicolumn{6}{c|}{Configuration C3 (left: $r=0.2$, right: $r=0.8$)}\\ \hline
    \textbf{Median} &1 & 1 & 1.031 &1 & 1&1.09  \\\hline
    \textbf{Mean} &2.06 &2.017 & 1.009 &1.983 & 1.927&1.059 \\\hline
    \textbf{Huber} &1.131 &  1.138& 1.037 & 1.102& 1.106& 1.095 \\  \hline
     \textbf{Lp-Median; $p=3$} & 4.22E+3& 4.24E+3  &1.034&4.54E+3&4.55E+3 & 1.092 \\    \hline
   \textbf{Lp-Median; $p=4$} & 3.56E+5  & 3.47E+5 & 1.006& 2.32E+5&2.28E+5  & 1.075 \\    \hline
   \textbf{Lp-Median; $p=5$} &  1.37E+6 & 1.33E+6  & 1.000 &  2.92E+6&2.81E+6    & 1.049 \\    \hhline{|=======|}
\multicolumn{1}{|l||}{} &  \multicolumn{6}{c|}{Configuration C4 (left: $r=0.2$, right: $r=0.8$)}\\ \hline
    \textbf{Median} &1& 1& 1.042  &1 &  1&1.099\\\hline
    \textbf{Mean} & 2.079 &2.034 & 1.019 &2.1 &2.042& 1.069 \\\hline
    \textbf{Huber} &1.136 & 1.44& 1.049 &1.174 & 1.178& 1.103\\\hline
     \textbf{Lp-Median; $p=3$} & 3.25E+3 & 3.28E+3& 1.052& 3.03E+3 &3.04E+3 & 1.101 \\    \hline
   \textbf{Lp-Median; $p=4$} &  1.40E+5 &1.42E+5 & 1.051 & 9.48E+4 &9.44E+4& 1.094 \\    \hline
   \textbf{Lp-Median; $p=5$} & 2.27E+6&2.29E+6& 1.051  & 9.86E+5 &9.75E+5  &   1.086  \\    \hline
    \end{tabular}}%
  \label{tab3}%
\end{table}%

\begin{table}[t!]
  \centering
  \caption{Cauchy distribution: Relative efficiencies w.r.t. to the (weighted) median and $\widehat{E}_f=(\frac{\text{det}(\widehat{\Sigma}_{\tht})}{\text{det}(\widehat{\Sigma}_{\tht}^w)})^{\frac{1}{2}}$}
    \renewcommand{\arraystretch}{1.2}
\scalebox{0.82}{   \begin{tabular}{|l||c|c|c||c|c|c|}
\cline{2-7}
 \multicolumn{1}{c|}{}&  \multicolumn{1}{c|}{$\mathbf{(\frac{\text{\bf det}(\widehat{\Sigma}_{\tht}^w) }{\text{\bf det}(\widehat{\Sigma}_{\text{\bf opt}}^w)})^{\frac{1}{2}}}$} & \multicolumn{1}{c|}{$\mathbf{(\frac{\text{\bf det}(\widehat{\Sigma}_{\tht})}{\text{\bf det}(\widehat{\Sigma}_{\text{\bf opt}})})^{\frac{1}{2}}}$} &  \multicolumn{1}{c|}{$\mathbf{\widehat{E}_f}$}  &  \multicolumn{1}{c|}{$\mathbf{(\frac{\text{\bf det}(\widehat{\Sigma}_{\tht}^w) }{\text{\bf det}(\widehat{\Sigma}_{\text{\bf opt}}^w)})^{\frac{1}{2}}}$} & \multicolumn{1}{c|}{$\mathbf{(\frac{\text{\bf det}(\widehat{\Sigma}_{\tht})}{\text{\bf det}(\widehat{\Sigma}_{\text{\bf opt}})})^{\frac{1}{2}}}$} &  \multicolumn{1}{c|}{$\mathbf{\widehat{E}_f}$}\\  \hhline{~|======|}
\multicolumn{1}{l|}{}&\multicolumn{6}{c|}{Configuration C1 (left: $r=0.2$, right: $r=0.8$)}\\\hline
\textbf{Median} &1 &1 & 2.246&1 & 1 & 3.868 \\\hline
\textbf{Mean} &2.98E+5 &2.13E+5 &1.605 &2.50E+5 & 1.42E+5   &2.191  \\\hline
   \textbf{Huber} & 1.161 &1.203 & 2.327  &1.181 &1.2& 3.928\\ \hline
 \textbf{Lp-Median; $p=3$} &8.86E+10 &4.72E+10 & 1.197 &8.07E+10 &3.27E+10 &1.570\\    \hline
   \textbf{Lp-Median; $p=4$} &6.13E+11 &3.18E+11 & 1.166 &3.404E+11 &1.19E+11  & 1.353   \\    \hline
   \textbf{Lp-Median; $p=5$} & 2.41E+11 &1.20E+11  & 1.121  & 8.59E+10& 2.76E+10    & 1.244 \\
    \hhline{|=======|}
\multicolumn{1}{|l||}{} &  \multicolumn{6}{c|}{Configuration C2 (left: $r=0.2$, right: $r=0.8$)}\\ \hline
    \textbf{Median} &1 &1 & 1.107 & 1 & 1 &1.309  \\ \hline
  \textbf{Mean} &7.99E+6  & 8.133+6   &1.126 &8.62E+6 &  8.16E+6  & 1.239 \\\hline
    \textbf{Huber} &1.161 & 1.170 & 1.115  &1.169 & 1.176 & 1.317 \\ \hline
     \textbf{Lp-Median; $p=3$} & 2.06E+12 &1.946E+12    & 1.043 &1.33E+12& 1.278E+12 &1.252 \\    \hline
   \textbf{Lp-Median; $p=4$} & 7.90E+10 &7.55E+10 & 1.058 &4.48E+10 &3.48E+10  & 1.014 \\    \hline
   \textbf{Lp-Median; $p=5$} &1.16E+11 &1.06E+11     & 1.009  & 5.66E+9  &6.00E+9 & 1.387  \\    \hhline{|=======|}
\multicolumn{1}{|l||}{}&  \multicolumn{6}{c|}{Configuration C3 (left: $r=0.2$, right: $r=0.8$)}\\\hline
    \textbf{Median} &1 & 1 &  1.030 &1 & 1 & 1.089 \\\hline
      \textbf{Mean} & 8.27E+5  &8.03E+5  & 1.000 &7.98E+5 &7.90E+5     &1.079 \\\hline
    \textbf{Huber} & 1.153& 1.156&  1.033& 1.152 &1.154 &1.092\\ \hline
     \textbf{Lp-Median; $p=3$} & 8.62E+11&8.46E+11    & 1.011& 2.65E+11 &2.56E+11 & 1.052  \\    \hline
   \textbf{Lp-Median; $p=4$} &2.82E+12 &2.74E+12 &  1.000 &4.92E+11 &4.70E+11   & 1.041  \\    \hline
   \textbf{Lp-Median; $p=5$} & 7.44E+11 &7.24E+11    &  1.002  & 7.96E+10 & 7.56E+10  & 1.034 \\   \hhline{|=======|}
 \multicolumn{1}{|l||}{} &   \multicolumn{6}{c|}{Configuration C4 (left: $r=0.2$, right: $r=0.8$)}\\\hline
    \textbf{Median} &1 & 1 & 1.04 & 1 & 1& 1.098 \\\hline
      \textbf{Mean} & 1.55E+6 &1.55E+6 &1.035 &1.39E+6 &1.38E+6     &1.087 \\\hline
    \textbf{Huber} &  1.129& 1.133&   1.043 &1.167 &1.169 & 1.100 \\\hline
     \textbf{Lp-Median; $p=3$} & 3.75E+12 &3.73E+12   & 1.034&  2.63E+12 &2.58E+12 & 1.080 \\    \hline
   \textbf{Lp-Median; $p=4$} &  3.74E+13& 3.70E+13&  1.027 & 1.86E+13 &1.81E+13 &  1.070 \\    \hline
   \textbf{Lp-Median; $p=5$} &  3.88E+13 &3.80E+13    & 1.020   &  1.39E+13 &1.34E+13    &1.059  \\    \hline
 \end{tabular}}%
  \label{tab4}%
\end{table}

In tables \ref{tab2} to \ref{tab4}, we report three different measures of efficiency: $(\frac{\text{det} (\widehat{\Sigma}_{\tht}^w)}{\text{det} (\widehat{\Sigma}_{\text{opt}}^w)})^{\frac{1}{2}}$ (weighted case), $(\frac{\text{det}(\widehat{\Sigma}_{\tht})}{\text{det}(\widehat{\Sigma}_{\text{opt}})})^{\frac{1}{2}}$ (unweighted case) and, $E_f$ defined in \eqref{e35}, estimated by $\widehat{E}_f= (\frac{\text{det}(\widehat{\Sigma}_{\tht})}{\text{det}(\widehat{\Sigma}_{\tht}^w)})^{\frac{1}{2}}$.  For the Gaussian distribution, $\widehat{\Sigma}_{\text{opt}}$ and $\widehat{\Sigma}_{\text{opt}}^{w}$ refer to the estimated variance of the (weighted) mean while for Student distributions, the (weighted) median is the reference. Not surprisingly, Huber estimators appear  more robust toward the distribution. Also, we note that, in the case of a bivariate Gaussian distribution with $r= 0.2$ and $r=0.8$  (Table \ref{tab2}),  the relative efficiency of all optimally weighted M-estimators is improved (compared to their unweighted version) whatever the configuration of clusters. However,   the quality improvement depends both on  configurations of clusters and  the value of $\rho$. For example, for the first configuration C1 (with one very big cluster), the efficiency of weighted estimators is doubled with respect to cases C2, C3 and C4. Namely, in the case $r=0.2$, the relative efficiency of the weighted  Huber estimator decreases from  2.461 for C1 to  1.038 for  C3. This highlights the impact of clusters' sizes on the variance: the improvement is even better when sizes are heterogeneous and big clusters are present.  The same behaviour is  observed for all tested estimators and they exhibit similar values for $\widehat{E}_f$. Next, in presence of a strong correlation ($r=0.8$), the relative efficiency is also improved: weighted M-estimators have a smaller variance than their unweighted version.  Again, presence of big clusters emphasizes this phenomenon. For example, the efficiency of the weighted Huber estimator decreases from  4.019 for C1 to  1.095 for  C3. Finally,  we get similar results for the Cauchy and the Student distributions (see tables \ref{tab3} and \ref{tab4}), where for ease of comparison, the same estimators had been computed even if some of them do not converge.

\section{Breakdown point}
\subsection{Computation of the breakdown point}
Results of the previous section show that one may  adapt weights to minimize the variance of the M-estimators. The benefit over the unweighted version is especially important for big clusters with high intra-correlation where the best relative efficiency is achieved by underweighting them. A consequence is that potential outliers in these clusters  have  reduced impact.  A possible measure of such robustness is the breakdown point. We recall here the definition given in \citet{DH83} (see also \citet{DG05b} for a discussion paper around this notion).
\begin{defi}
\label{d41}
The finite sample replacement breakdown point of  $\widehat{\tht}_n(X)$ built with  $n$ observations is defined by:
$$
\epsilon^{\ast}_n=\min_{1\leq k\leq n}\left\{\frac{k}{n} : \sup_{Y_k}{\left\|\widehat{\theta}_n(Y_k)-\widehat{\tht}_n(X)\right\|=\infty}\right\}
$$
where $Y_k$  denotes the corrupted sample from $X$, obtained by replacing $k$ points of $X$ with arbitrary values.
\end{defi}
We denote by  $\widehat{\tht}_n(X)$ and $\widehat{\tht}_n^w(X)$ the unweighted and weighted M-estimators based on $X$. We reorganize the indexation in
 $X=\{X_1,\dotsc,X_{N_n}\}$  with associated weights $\{w_1,\dotsc,w_{N_n}\}$, so these estimators can be written as
\begin{equation*}
\widehat{\tht}_n(X)=\argmin{a\in\Theta} \frac{1}{N_n}\sum^{N_n}_{i=1}{\rho(X_{i},a)}
\text{ and }
\widehat{\tht}_n^w(X)=\argmin{a\in\Theta} \frac{1}{N_n}\sum^{N_n}_{i=1}{w_{i} \rho(X_{i},a)}.
\end{equation*}
Their breakdown points are $\epsilon^{\ast}_{N_n}=\min\limits_{1\leq k\leq N_n}\left\{\frac{k}{N_n} : \sup\limits_{Y_k}{\left\|\widehat{\theta}_n(Y_k)-\widehat{\tht}_n(X)\right\|=\infty}\right\} $  and
\begin{equation}
\label{e51}
\epsilon^{w\ast}_{N_n}=\min_{1\leq k\leq N_n}\left\{\frac{k}{N_n} : \sup_{Y_k}{\left\|\widehat{\theta}^w(Y_k)-\widehat{\tht}_n^w(X)\right\|=\infty}\right\}
\end{equation}
where $Y_k$ is again the corrupted sample from $X$ with $k$ arbitrary values. In this finite framework, we suppose that $\frac{1}{N_n}\sum_{i=1}^n\sum_{j=1}^{m_i}w_{ij}=1$, that can be written as  $\sum_{i=1}^{N_n} w_i = N_n$ with the new indexation. The following result gives two bounds for $\epsilon^{w\ast}_{N_n}$ depending on the weights. For conciseness and clarity, we choose them as  rational: $w_i= \frac{\ell_i}{L}$ for some $L\ge 1$.
\begin{theo}
\label{t42}
Suppose that the breakdown point of  $\widehat{\tht}_n(X)$ is such that:  for all integers $m\ge 1$,  $\epsilon^\ast_{N_n} \le \epsilon^\ast_{mN_n} \le \epsilon_0^\ast$  where $\epsilon_0^\ast$ represents the asymptotic breakdown point. The estimator $\widehat{\tht}_n^w(X)$ has a breakdown point $\frac{k_w^\ast}{N_n}$ with $k_w^\ast\in [[k_{\sste N}^{\ast}, k_0^\ast]]$  and the bounds  $k_{\sste N}^{\ast}$ and $k_0^\ast$ are respectively defined by:
\begin{align*}
k_{\sste N}^{\ast}&=\min_{1\leq k\leq N_n}\left\{k :  \ex  i_1,\dotsc,i_k \in \{1,\dotsc,N_n\}~\Big|~w_{i_1}+\dotsb+w_{i_k}\geq \epsilon^{\ast}_{N_n} N_n \right\},\\
k_0^\ast&=\min_{1\leq k\leq N_n}\left\{k : \ex i_1,\dotsc,i_k \in \{1,\dotsc,N_n\}~\Big|~w_{i_1}+\dotsb+w_{i_k}\geq \epsilon^{\ast}_{0} N_n \right\}.
\end{align*}
\end{theo}
\begin{proof}
We follow the main steps of the proof given in \citet{NLO06} for the breakdown point of the spatial median. Clearly $k$ corresponds to the number of variables in $Y_k$ not present in $X$: $ k= \#\{Y_k\setminus(Y_k\cap X) \}$ with $\#\{A\}$  the cardinal of $A$. From the definition of $\widehat{\tht}_n(X)$,  we get that  $
\sup\limits_{Y_k}{\left\|\widehat{\theta}_n(Y_k)-\widehat{\tht}_n(X)\right\|=\infty}$ is equivalent to $ \{k \geq \epsilon^{\ast}_{N_n} N_n\}$, so that:
\begin{equation}\label{e52}
\sup_{Y_k}{\left\|\widehat{\theta}_n(Y_k)-\widehat{\tht}_n(X)\right\|=\infty} \Longleftrightarrow \big\{\#\{Y_k\setminus(Y_k\cap X) \} \geq \epsilon_{N_n}^{\ast} \#\{X\}\big\}.
\end{equation}
For $w_i=\frac{\ell_i}{L}$,   $i=1,\dotsc,N_n$ and $\ell_i,L \in \mathds{N}^\ast$, the weighted M-estimator  can be written as $
\widehat{\tht}_n^w(X)=\argmin{a\in\Theta} \frac{1}{N_nL}\sum^{N_n}_{i=1}{\ell_i \rho(X_{i},a)}$. Therefore $\widehat{\tht}_n^w(X)$, associated with $X$, is also the unweighted estimator $\widehat{\tht}_n(\widetilde{X})$ where $\widetilde{X}$
is defined by each  $X_i$ of $X$ repeated  $\ell_i$ times (and similarly  for the set $\widetilde{Y}_k$  deduced from $Y_k$). These transformations allow us to write the breakdown point given by \eqref{e51} as:
\begin{equation*}
\epsilon^{w\ast}_{N_n} =\min_{1\leq k\leq N_n}\left\{\frac{k}{N_n} : \sup_{\widetilde{Y}_k}{ \left\| \widehat{\tht}_n (\widetilde{Y}_k) - \widehat{\tht}_n (\widetilde{X})\right\| =\infty} \right\}.
\end{equation*}
Next, using \eqref{e52} and the condition   $\sum_{i=1}^{N_n}\ell_i=LN_n$, we obtain
\begin{equation*}
\epsilon^{w\ast}_{N_n} =\min_{1\leq k\leq N_n}\left\{\frac{k}{N_n} : \#\{\widetilde{Y}_k \setminus (\widetilde{Y}_k\cap \widetilde{X}) \} \geq \epsilon_{LN_n}^{\ast} \#\{\widetilde{X}\} \right\},
\end{equation*}
where, if $X_{i_1},\dotsc,X_{i_k}$ are the $k$ points replaced in $X$, one has to replace $\ell_{i_1}+\dotsb+\ell_{i_k}$ points in  $\widetilde{X}$ by arbitrary values to obtain $\widetilde{Y}_k$. Moreover $\#\{\widetilde{X}\}=\sum_{i=1}^{N_n}\ell_i=LN_n$,  so the breakdown point is given by
\begin{align} \notag
\epsilon^{w\ast}_{N_n} &=\min_{1\leq k\leq N_n}\left\{\frac{k}{N_n} : \ex i_1,\dotsc,i_k \in \{1,\dotsc,N_n\}~\Big|~\ell_{i_1}+\dotsb+\ell_{i_k}  \geq \epsilon_{LN_n}^{\ast} \sum_{i=1}^{N_n} \ell_i \right\}\\
\label{e53}
&=\min_{1\leq k\leq N_n}\left\{\frac{k}{N_n} : \ex i_1,\dotsc,i_k \in \{1,\dotsc,N_n\}~\Big|~w_{i_1}+\dotsb+w_{i_k}  \geq \epsilon_{LN_n}^{\ast} N_n \right\}.
\end{align}
Let $k^{\ast}_w$ be the minimal value obtained in \eqref{e53}. By definition of $k_0^{\ast}$,
 the property  $\epsilon^{\ast}_{0} \ge  \epsilon_{LN_n}^{\ast}$ implies
that $w_{i_1}+\dotsb + w_{i_{k_0^{\ast}}} \ge  \epsilon^{\ast}_{LN_n} N_n$. Therefore, we may deduce that $k^{\ast}_w \le k_0^{\ast}$. In the same way, from the definition of $k_{\sste N}^{\ast}$ and the condition  $\epsilon^{\ast}_{N_n} \le  \epsilon_{LN_n}^{\ast}$, we also get  $k^{\ast}_w \ge k_{\sste N}^{\ast}$.\end{proof}
\begin{rmk}  For the spatial median, one has $\ve_{N_n}^{\ast}=\frac{\left\lfloor \frac{N_n-1}{2}\right\rfloor}{N_n}$ \citep{CR92}. In this case, $\epsilon_0^\ast= \frac{1}{2}$  and the condition $\epsilon^\ast_{N_n} \le \epsilon^\ast_{mN_n}\le \epsilon_0^\ast$ is satisfied for all $m\ge 1$. In dimension $d$ with  $\ve_{N_n}^{\ast}=\frac{\left\lfloor \frac{N_n-d+1}{2}\right\rfloor}{N_n}$, the condition is fulfilled too. \end{rmk}

From Theorem~\ref{t42}, the breakdown point of a weighted M-estimator depends more on its weights than on potential outliers. Furthermore, the proof shows that for $w_i= \frac{\ell_i}{L}$, its exact expression takes the form given in \eqref{e53}.
It is worth noting that if  $\epsilon_{LN_n}^{\ast}$ and $\epsilon_0^{\ast}$ are very closed, we get
$$\epsilon^{w\ast}_{N_n}=\min_{1\leq k\leq N_n}\left\{\frac{k}{N_n} : \ex  i_1,\dotsc,i_k \in \{1,\dotsc,N_n\}~\Big|~w_{i_1}+\dotsb+w_{i_k}  \geq \epsilon_{0}^{\ast} N_n \right\}$$
which generalizes the definition given by \citet{NLO06} for the weighted spatial median (as they use the asymptotic breakdown point of 0.5 to derive their results). We conclude  with two remarks: weights do not improve the asymptotic breakdown point of the unweighted case, and, a trade off   between optimal efficiency and maximal breakdown point, should be found.

\begin{rmk} Suppose that the weighted estimator achieves its maximal breakdown point, $\epsilon_{N_n}^{w\ast}= \frac{k_0^{\ast}}{N_n}$, with $k_0^{\ast}$ such that
$$
k_0^\ast=\min_{1\leq k\leq N_n}\left\{k : \ex i_1,\dotsc,i_k \in \{1,\dotsc,N_n\}~\Big|~  w_{i_1}+\dotsb+w_{i_k}\geq \epsilon^{\ast}_{0} N_n \right\}.
$$
If the weights $w_i$ are ranked in ascending order: $w_{(1)}<\dotsb <w_{(N_n)}$, the minimality of $k$ is guaranteed by  replacing the observations with the largest $k_0^{\ast}$ weights where $k_0^\ast$ satisfies
\begin{equation}\label{e54}
\begin{cases} w_{(N_n)}+\dotsb+w_{(N_n -k_0^{\ast}+1)} \ge \epsilon^{\ast}_{0} N_n \\
w_{(N_n)}+\dotsb+w_{(N_n -k_0^{\ast}+2)} < \epsilon^{\ast}_{0} N_n.
\end{cases}
\end{equation}
\end{rmk}
From the last remark, the minimal improvement of the unweighted  breakdown point corresponds to $k_0^{\ast}= \epsilon_0^{\ast}N_n+1$ where, for the sake of clarity, we choose  $N_n$  such that $\epsilon_0^{\ast}N_n$ is again an integer. The second part of  \eqref{e54} becomes $\sum_{i=0}^{\epsilon_0^{\ast}N_n -1} w_{(N_n-i)} < \epsilon_0^{\ast}{N_n}$.  As this sum includes $\epsilon_0^{\ast}N_n$  terms, necessarily one gets that $w_{(N_n - \epsilon_0^{\ast}N_n+1)} <1$, so $w_{(1)}\le \dotsb \le w_{(N_n - \epsilon_0^{\ast}N_n)}\le w_{(N_n - \epsilon_0^{\ast}N_n+1)} < 1$. This leads us to a contradiction: on one hand $\sum_{i=1}^{N_n} w_{(i)}=N_n$, and on the other hand, $\sum_{i=1}^{N_n - \epsilon_0^{\ast}N_n} w_{(i)} <  N_n - \epsilon_0^{\ast}N_n$ and $\sum_{i=N_n - \epsilon_0^{\ast}N_n+1}^{N_n} w_{(i)}=\sum_{i=0}^{\epsilon_0^{\ast}N_n -1} w_{(N_n-i)} <  \epsilon^{\ast}_{0} N_n $. The asymptotic breakdown point of the unweighted case could not be increased.
\begin{rmk}
If the maximal breakdown point of the unweighted version (obtained with $w_i\equiv1$) is a proportion $\epsilon_0^{\ast}$ (with $0< \epsilon_0^{\ast}\le \frac{1}{2}$) of the data, then from  the first equation of \eqref{e54}, we have to overweight (with weights not smaller than one) $\epsilon_0^{\ast}N_n$ observations  in $\widehat{\tht}_n^w(X)$ to reach  this value. \end{rmk}
We conclude that,   at least in the case where $w_{ij}\equiv w_i$, there is  no hope to simultaneously maximize the  breakdown point and the relative efficiency (since the smallest weights are assigned to the biggest clusters and, consequently to the maximal number of variables). This fact is illustrated in the following section.
\subsection{Numerical results}

\begin{table}[b!]
  \centering
  \caption{Breakdown point $\epsilon^{w\ast}_n$ for the weighted spatial median and the weighted Huber estimator}
    \begin{tabular}{|c|c|rrrr|}
    \cline{3-6}
 \multicolumn{2}{c|}{} & \multicolumn{4}{c|}{Weighted spatial median}  \\
\cline{3-6}
 \multicolumn{2}{c|}{}    &C1     & C2     & C3     & C4\\ \hline
\multirow{2}{*}{$r= 0.2$} &
    $w_{(N_n)}+\ldots+w_{(N_n-k-1)}$ & 51.921 & 50.173 & 50.817 & 50.081 \\
    \cline{2-6}
  &  $ \epsilon^{w\ast}_{N_n}$ & 23\%  & 37\%  & 46\%  & 41\% \\ \hline
\multirow{2}{*}{$r= 0.8$} &
    $w_{(N_n)}+\ldots+w_{(N_n-k-1)}$ & 52.099 & 50.079 & 50.177 & 50.293 \\
    \cline{2-6}
  &  $\epsilon^{w\ast}_{N_n}$  & 21\%  & 23\%  & 41\%  & 36\%   \\ \hline
   \multicolumn{6}{c}{} \\
      \cline{3-6}
 \multicolumn{2}{c|}{} & \multicolumn{4}{c|}{Weighted Huber estimator} \\
\cline{3-6}
 \multicolumn{2}{c|}{} & C1     & C2     & C3     & C4 \\ \hline
\multirow{2}{*}{$r= 0.2$} &
    $w_{(N_n)}+\ldots+w_{(N_n-k-1)}$  &50.118 & 50.326 & 50.055 & 50.611 \\
    \cline{2-6}
  &  $ \epsilon^{w\ast}_{N_n}$ & 22\%     & 36\%      & 45\%      & 41\% \\ \hline
\multirow{2}{*}{$r= 0.8$} &
    $w_{(N_n)}+\ldots+w_{(N_n-k-1)}$ & 51.985 & 50.443 & 50.117 & 50.219 \\
    \cline{2-6}
  &  $\epsilon^{w\ast}_{N_n}$  & 21\%      & 23\%      & 41\%      & 36\%   \\ \hline
    \end{tabular}%
  \label{tab7}%
\end{table}%

In this part we evaluate the breakdown point of the weighted spatial median and the weighted Huber estimator, whose unweighted versions have a maximal breakdown point of $0.5$. We consider the configurations C1-C4 defined in section~\ref{subsec:efficiency}
and  a centered bivariate Gaussian distribution with  $r=0.2$ and $r=0.8$. We select the optimal weights (with $w_{ij}\equiv w_i$, $i=1,\dotsc,10$) maximizing the relative efficiency (see Table \ref{tab1}).
The  breakdown points, computed for these two estimators, are presented in Table~\ref{tab7}: in each case, they are far less than  $50\%$. Variations are observed according to the number of variables by clusters and  value of correlation.
Both a strong correlation and the presence of big clusters worsen the breakdown point (see for example the case of configuration C1). The most favorable configuration seems to be C3: it corresponds to a maximal size of clusters limited to 12. We conclude that optimal weights improve significantly the efficiency but can drastically reduce the breakdown point.

\appendix
\renewcommand*{\thesection}{\Alph{section}}

\section{Proof of the auxiliary results}

\subsection{Proof of Lemma~\ref{l31}}
We make use of two limit theorems for independent random variables: the weak law of
\citet[p. 356]{CT97}  and the strong one of \citet{OS02}. For  the  sake  of  completeness,  we  state below their  main  results
in  a  version  suitable  for  our  purposes.

\begin{theo}\label{ta1}  For each $n \ge 1$, let ${Y_{in}, i=1,\dotsc,n}$ be
independent r.v.s and $S_n = \sum_{i=1}^n Y_{in}$.
\begin{enumerate}[(1)]
\item If $\sum\limits_{i=1}^n \esp(\abs{Y_{in}}^{1+\delta})\tv[n\to\infty]{} 0$ for some $\delta>0$, then $S_n -  \sum\limits_{i=1}^n \esp(Y_{in})$ converges to 0 in probability.
\item If $\esp(Y_{in})=0$, $\sum\limits_{n\ge 1} \sum\limits_{i=1}^n \esp(\abs{Y_{in}}^{1+\delta}) <\infty$ for some $\delta>1$,   and for some $s>0$,  $\sum\limits_{n\ge 1} \big(\sum\limits_{i=1}^n   \esp(\abs{Y_{in}}^{2})\big)^s <\infty$, then $S_n$ converges to 0 almost surely.
\end{enumerate}
\end{theo}
Remark that in the original version of \citeauthor{CT97}, the result is given in the form $S_n -  \sum_{i=1}^n \esp(Y_{in}\indi_{\{\abs{Y_{in}} <1\}})\to 0$ in probability. Applying H\"{o}lder and Markov inequalities, one obtains easily that $ \sum_{i=1}^n \esp(Y_{in}\indi_{\{\abs{Y_{in}} <1\}})= \sum_{i=1}^n \esp(Y_{in}) + o(1)$ if  $\sum_{i=1}^n \esp(\abs{Y_{in}}^{1+\delta})\to 0$ as $n\to\infty$.

We begin to establish  the consistency of $M_n^w(a)$ to $M(a)=\esp(\rho(X_{11},a))$ for all $a\in\Theta$. To this end, we set $Y_{in}= \frac{1}{N_n}\sum_{j=1}^{m_i} w_{ij} \rho(X_{ij},a)$, these random variables are independent with expectation $\frac{1}{N_n}\sum_{j=1}^{m_i}w_{ij}\espt(\rho(X_{11},a))$. Using Minkowski inequality, we obtain:
\begin{equation}
\esp(Y_{in}^{1+\delta}) \le \big(\sum_{j=1}^{m_i} \frac{w_{ij}}{N_n}\big)^{1+\de} \sup_{a\in\Theta} \esp(\rho^{1+\de}(X_{11},a)). \tag{a.1} \label{ea1}
\end{equation}
From conditions A\ref{a31}(c)-(i), we may deduce that $\sum\limits_{i=1}^n \esp(\abs{Y_{in}}^{1+\delta})\tv[n\to\infty]{} 0$.
Then, for all  $a\in \Theta$, the first part of Theorem~\ref{ta1} gives:
\begin{equation} \tag{a.2}\label{ea2}
\frac{1}{N_n}\sum^n_{i=1} \Big(\sum^{m_i}_{j=1}{w_{ij} \rho(X_{ij},a)} -  \sum_{j=1}^{m_i} w_{ij}\esp\big(\rho(X_{11},a) \big)\Big)\xrightarrow[n\to\infty]{p}0.
\end{equation}
By A\ref{a31}(b), we have $\frac{1}{N_n} \sum\limits_{i=1}^n \sum\limits_{j=1}^{m_i} w_{ij}\tv[n\to\infty]{} 1$, so the result: $M_n^w(a)\tv[n\to\infty]{p} M(a)$.
\\
For the strong consistency, \eqref{ea1} with $\delta\ge 1$  shows  that $\sum\limits_{n\ge 1} \sum\limits_{i=1}^n \esp(\abs{Y_{in}-\esp(Y_{in})}^{1+\delta})$   and   $\sum\limits_{n\ge 1} \big(\sum\limits_{i=1}^n   \esp(\abs{Y_{in}-\esp(Y_{in})}^{2})\big)^s$ are convergent series under the condition A\ref{a31}(c)-(ii).  From Theorem~\ref{ta1}-(2), we obtain that:
\begin{equation}\tag{a.3}\label{ea3}
\frac{1}{N_n}\sum^n_{i=1} \Big(\sum^{m_i}_{j=1}{w_{ij} \rho(X_{ij},a)} -  \sum_{j=1}^{m_i} w_{ij}\esp\big(\rho(X_{11},a) \big)\Big)\xrightarrow[n\to\infty]{a.s.}0.
\end{equation}
and the convergence of $\sum\limits_{i=1}^n\esp(Y_{in})$ to $M(a)$ yields the final result.

Now, we turn to the uniform convergence. From the compactness of $\overline{\Theta}$ (where $\overline{A}$ denotes the closure of the set $A$), we get that $\Theta\subset\bigcup^{r_n}_{l=1}B(a_{ln},h_n)$, with $B(a_{ln},h_n)$ the open ball of center  $a_{ln}$ and radius $h_n\to 0$. Then,
\begin{multline}\tag{a.4}\label{ea4}
\sup_{a\in\Theta} \left|M_n^w(a)-M(a)\right| \le \max\limits_{l=1,\dots,r_n} \sup_{a\in B(a_{ln},h_n)} \left|M_n^w(a)-M(a)\right|  \\\le \max\limits_{l=1,\dots,r_n} \Big(\left|M_n^w(a_{ln})-M(a_{ln})\right| + \sup_{a\in B(a_{ln},h_n)}  \left|M_n^w(a)-M_n^w(a_{ln})\right| \\+  \left|M(a_{ln})-M(a)\right|\Big).
\end{multline}
Since $\rho(x,\cdot)$ is $k(x)$-H\"{o}lderian, we get the upper bound
$$\abs{M_n^w(a)-M_n^w(a_{ln})} \le \frac{1}{N_n}\sum^n_{i=1}\sum^{m_i}_{j=1}{w_{ij} k(X_{ij})} \norm{a-a_{ln}}^{\lambda}, \; \lambda\in]0,1].$$
Next from the Theorem~\ref{ta1}-(1) and conditions A\ref{a31}(b),(c)-(ii), we obtain that $$\frac{1}{N_n}\sum^n_{i=1}\sum^{m_i}_{j=1}{w_{ij} k(X_{ij})}\xrightarrow[n\to\infty]{p}\espt \big(k(X_{11})\big).$$
We may deduce that  for $n$ large enough, the first  term of \eqref{ea4} is a ${\cal O}(h_n^{\lambda})$ uniformly in $l$.  The almost sure consistency follows in a similar way from the condition A\ref{a31}(c)-(ii). Moreover, the H\"{o}lderian condition on $M(a)$ implies the same result for $\abs{M(a_{ln})-M(a)}$.  Finally, we study  the convergence of the first  term $\abs{M_n^w(a_{ln})-M(a_{ln})}$. Now, setting $Y_{in} = \frac{1}{N_n}\sum_{j=1}^{m_i} w_{ij} \rho(X_{ij},a_{nl})$,  we  may follow the steps of the first part of the proof because the bound obtained in \eqref{ea1} is uniform in $a$. From \eqref{ea2}-\eqref{ea3} with $a$ replaced by $a_{nl}$, we obtain that
 $$
 \frac{1}{N_n}\sum_{i=1}^n \sum_{j=1}^{m_i} w_{ij} \rho(X_{ij}, a_{nl}) - \frac{1}{N_n} \sum_{i=1}^n \sum_{j=1}^{m_i} w_{ij} \espt \rho(X_{ij}, a_{nl})\tv[n\to\infty]{} 0$$
either in probability or almost surely. We conclude with
$$\abs{\esp(M_n^w(a_{nl})) - M(a_{nl})} \le \big|\frac{1}{N_n}\sum_{i=1}^n \sum_{j=1}^{m_i} w_{ij} - 1\big| \sup_{a\in\Theta} \esp(\rho(X_{11},a))=o(1)$$ uniformly in $l$. Collecting all the results, the uniform  convergence of $M_n^w$ to $M(a)$ follows from \eqref{ea4}.~

\subsection{Proof of Lemma~\ref{l32}}
We use  the multivariate version of the Lindeberg-Feller theorem, see  \citet[p.~147]{Ra73}. First for $i=1,\dotsc,n$, the vectors ${\xi_i^w:=\sum\limits^{m_i}_{j=1}{w_{ij} \psi(X_{ij},\theta)}}$  are independent, centred and with variances $V_{\tht,i}$, $i=1,\dotsc,n$,  defined by
$$V_{\tht,i}= \sum^{m_i}_{j=1}w_{ij}^2 \espt\psi(X_{ij},\theta)\psi^{\sste T}(X_{ij},\theta)
+\sum^{m_i}_{j\neq j'} w_{ij}w_{ij'}\espt\psi(X_{ij},\theta)\psi^{\sste T}(X_{ij'},\theta).$$
Since all pairs $(X_{ij'},X_{ij})$, $j \neq j'$,  have the same correlation in cluster $i$, we set $C_{\tht,i}=\espt\psi(X_{ij},\theta)\psi^{\sste T}(X_{ij'},\theta)$ and $B_{\tht}:=\espt\psi(X_{1j},\theta)\psi^{\sste T}(X_{1j},\theta)$.  Then the conditions A\ref{a32}(a)-(b) give:
$$
\frac{1}{N_n}\sum^n_{i=1}V_{\tht,i} \xrightarrow[n\to\infty]{~} c_{w}B_{\tht}+C_{\tht}^w.
$$
Next as $\lim_{n \to \infty} \frac{N_n}{n}=\ell$,  one can write:
$$
\frac{1}{n}\sum^n_{i=1}V_{\tht,i} \xrightarrow[n\to\infty]{~} \ell(c_{w}B_{\tht}+C_{\tht}^w):= \Sigma.
$$
Then, we verify the Lindeberg condition:
$$E_n:=\frac{1}{n}\sum^n_{i=1} \espt\big(\norm{\xi_i^w}^2 \indi_{\{\left\|\xi_i^w\right\|>\epsilon\sqrt{n}\}}\big) \tv[n\to\infty]{} 0.$$
Using the  H\"{o}lder, Minkowski and Markov inequalities, we may bound $E_n$ by
\begin{align*}
E_n\leq \espt\big(\left\|\psi(X_{11},\theta)\right\|^{2+\eta}\big)  \frac{1}{\epsilon^{\eta}n^{1+\frac{\eta}{2}}} \sum^n_{i=1} (\sum_{j=1}^{m_i} w_{ij})^{2+\eta}\xrightarrow[n\to\infty]{} 0
\end{align*}
with  the help of condition A\ref{a32}-(c) for all $\epsilon>0$. Then $\frac{1}{\sqrt{n}}\sum^n_{i=1}\xi_i^w\stackrel{d}{\rightarrow}N\left(0,l(c_{w}B_{\tht}+C_{\tht}^w)\right)$ so $\sqrt{N_n}T_n^w(\theta)$ is asymptotically normal with covariance matrix $c_{w}B_{\tht}+C_{\tht}^w$.

\hfill

\subsection{Proof of Lemma~\ref{l33}}
We have $\dot{T}_n^w(\theta)=\frac{1}{N_n}\sum^n_{i=1}\sum^{m_i}_{j=1}{w_{ij} \dot{\psi}(X_{ij},\theta)}$ and for some $\delta >0$,
$\espt(\dot{\psi}(X_{11},\tht))^{1+\de}<\infty$, so the  $P$-consistency of $\dot{T}_n^w(\theta)$ follows from the first part of Theorem~\ref{ta1} and the condition A\ref{a32}(d).

\subsection{Proof of Proposition~\ref{p41}}
Let us denote $\widehat{B}_n^{w(k,l)}(a)$, the $(k,l)$ element of the matrix $\widehat{B}_n^w(a)$,
$$\widehat{B}_n^{w(k,l)}(a)=\frac{1}{N_n} \sum_{i=1}^n \sum_{j=1}^{m_i}w_{ij}^2  \psi_k(X_{ij},a)\psi_l(X_{ij},a).
$$
Again we apply the first part of Theorem \ref{ta1} with $$Y_{in}^{(k,l)} =  \frac{1}{N_n}\sum_{j=1}^{m_i} w_{ij}^2 \psi_k(X_{ij},a)\psi_l(X_{ij},a).$$
With Minkowski and Cauchy Schwarz inequalities, we obtain the bound
$$
\sum_{i=1}^n \esp\abs{Y_{in}^{(k,l)}}^{1+\frac{\eta}{2}} \le(\esp\abs{\psi_k^{2+\eta}(X_{11},a)})^{\frac{1}{2}} \esp(\abs{\psi_l^{2+\eta}(X_{11},a)})^{\frac{1}{2}}\sum_{i=1}^n\big(\frac{\sum_{j=1}^{m_i}w_{ij}^2}{N_n}\big)^{1+\frac{\eta}{2}}=o(1).
$$
Next $\sum_{i=1}^n\esp(Y_{in}^{(k,l)})\to c_w \esp(\psi_k(X_{11},a)\psi_l(X_{11},a))$, hence the result.

Proofs are similar for the terms  $\widehat{C}_n^w(a)$ and  $\widehat{V}_n(a)$. In particular, the $P$-consistency of $\widehat{V}_n(a)$ follows with the conditions A\ref{a32}-(d).

\bibliographystyle{elsarticle-harv}
\bibliography{bibMEP}
\end{document}